\documentclass[11pt]{article}
\usepackage{authblk}

\usepackage{geometry}
\usepackage{amsmath,amssymb,mathrsfs}
\usepackage{color}
\usepackage{graphicx}
\usepackage{subfigure}
\usepackage{soul}
\usepackage{bbding}
\usepackage{lineno,hyperref}
\usepackage{multirow, booktabs, stmaryrd}
\modulolinenumbers[5]

\newcommand\btau{\boldsymbol{\tau}}
\DeclareMathOperator*{\Loss}{Loss}

\def\bp{\mathbf{p}}
\def\bx{\mathbf{x}}
\def\bn{\mathbf{n}}
\def\bu{\mathbf{u}}

\def\bg{\mathbf{g}}
\def\bF{\mathbf{F}}
\def\bX{\mathbf{X}}
\def\grad{\nabla}
\def\div{\nabla\cdot}
\def\laplace{\Delta}
\def\pd{\partial}
\newcommand{\beq}{\begin{equation}}
\newcommand{\eeq}{\end{equation}}
\newcommand{\beqs}{\begin{eqnarray}}
\newcommand{\eeqs}{\end{eqnarray}}
\newcommand{\beqsn}{\begin{eqnarray*}}
\newcommand{\eeqsn}{\end{eqnarray*}}
\newcommand{\bary}{\begin{array}}
\newcommand{\eary}{\end{array}}

\newcommand{\blue}[1]{\textcolor{black}{#1}}
\newcommand{\red}[1]{\textcolor{black}{#1}}
\newcommand*{\dbblk}[1]{\llbracket #1 \rrbracket}

\title{An efficient neural-network and finite-difference hybrid method for elliptic interface problems with applications}

\author[1,4]{Wei-Fan Hu}
\author[2,4]{Te-Sheng Lin}
\author[3]{Yu-Hau Tseng}
\author[2]{Ming-Chih Lai}

\affil[1]{Department of Mathematics, National Central University, Taoyuan 32001, Taiwan}
\affil[2]{Department of Applied Mathematics, National Yang Ming Chiao Tung University, Hsinchu 30010, Taiwan}
\affil[3]{Department of Applied Mathematics, National University of Kaohsiung, Kaohsiung 81148, Taiwan}
\affil[4]{National Center for Theoretical Sciences, National Taiwan University, Taipei 10617, Taiwan}

\begin{document}

\maketitle

\begin{abstract}
A new and efficient neural-network and finite-difference hybrid method is developed for solving Poisson equation in a regular domain with jump discontinuities on embedded irregular interfaces. Since the solution has low regularity across the interface, when applying finite difference discretization to this problem, an additional treatment accounting for the jump discontinuities must be employed. Here, we aim to elevate such an extra effort to ease our implementation by machine learning methodology. The key idea is to decompose the solution into singular and regular parts. The neural network learning machinery incorporating the given jump conditions finds the singular solution, while the standard \red{five-point Laplacian discretization} is used to obtain the regular solution with associated boundary conditions. Regardless of the interface geometry, these two tasks only require supervised learning for function approximation and a fast direct solver for Poisson equation, making the hybrid method easy to implement and efficient. The two- and three-dimensional numerical results show that the present hybrid method preserves second-order accuracy for the solution and its derivatives, and it is comparable with the traditional immersed interface method in the literature. As an application, we solve the Stokes equations with singular forces to demonstrate the robustness of the present method.\\

\end{abstract}

%%%%%%%%%%%%%%%%%%%%%%%%%%%%%%%%%%%%%%%%%
\section{Introduction}
%%%%%%%%%%%%%%%%%%%%%%%%%%%%%%%%%%%%%%%%%

In this paper, we aim to solve a $d$-dimensional ($d = 2$ or 3) elliptic interface problem defined in a regular domain $\Omega \subset \mathbb{R}^d$, which is separated by an embedded interface $\Gamma$ such that the subdomains inside and outside the interface are denoted by $\Omega^-$ and $\Omega^+$, respectively. Along the interface $\Gamma$, there
exists jump discontinuities that the solution must be satisfied. With the associated boundary condition, the problem
takes the form
\begin{align}
\laplace u(\bx) = f(\bx), \quad &\bx\in \Omega^- \cup \Omega^+, \label{Eq:Poisson}\\
\dbblk{u(\bx)} = \gamma(\bx), \quad \dbblk{\partial_n u(\bx)} = \rho(\bx), \quad &\bx\in\Gamma, \label{Eq:jump}\\
u(\bx) = u_b(\bx), \quad &\bx\in\partial\Omega. \label{Eq:bdc}
\end{align}
Here, the jump $\dbblk{\cdot}$ indicates the quantity approaching from $\Omega^+$ side minus the one from $\Omega^-$ side; the shorthand $\partial_n u$ represents the normal derivative $\grad u\cdot\bn$ in which $\bn$ is the normal vector pointing from $\Omega^-$ to $\Omega^+$.
Notice that, here the underlying differential equation is subject to the Dirichlet-type boundary condition for illustration purpose, while other types of boundary condition (Neumann or Robin) will not change the main ingredients presented here. Since the Poisson equation is considered in Eq.~(\ref{Eq:Poisson}), we simply call the above problem as the Poisson interface problem hereafter.

As seen from Eqs.~(\ref{Eq:Poisson})-(\ref{Eq:bdc}), the solution and its partial derivatives have jumps across the interface.  So, when applying the finite difference discretization to this problem, an additional treatment accounting for those jump discontinuities must be employed at the grid points near the interface. Over the past few decades, different discretization methodologies have been successfully developed to capture those jump conditions sharply or to improve the overall numerical accuracy, such as the immersed interface method (IIM)~\cite{HHL19, HLY15, LT08, LL94}, ghost fluid method (GFM)~\cite{EG20, LFK20}, Voronoi interface method~\cite{GLTG15}, to name a few. \blue{Different approaches for solving interface problems such as the immersed finite element method (IFEM)~\cite{GLL08, HLL11}} or other methods can be found in~\cite{LI2006} and the references therein.

On the other hand, much attention has recently been paid to applying deep neural networks~(DNNs) to solve elliptic interface problems, rather than using traditional numerical methods to solve such problems. Despite the success of the two mainstream deep learning approaches (Physics-Informed Neural Networks (PINNs)~\cite{RPK19, SDK20} and the deep Ritz method~\cite{EY18}) in solving partial differential equations with smooth solutions, learning methods based on these two frameworks for solving elliptic interface problems with jump discontinuities remain to be improved. The main and intrinsic difficulty may be attributed to the fact that the usual activation functions used in DNNs are generally smooth; thus, DNN function approximators seem to be incapable of representing discontinuous functions. \blue{To approximate such discontinuous solutions (or functions) and tackle the elliptic interface problems, multiple independent networks need to be established and linked with each other by imposing the jump conditions, see, e.g., piecewise DNNs~\cite{HHM22}, interfaced neural networks~\cite{WL22}, and deep unfitted Nitsche method~\cite{GY22}.} The resulting prediction errors in their test examples reach the magnitude $O(10^{-3})$ to $O(10^{-4})$ in relative $L^2$ norm. Moreover, training these DNN models comes at the cost of having to train a separate neural network in each subdomain independently. Until very recently, the authors of this paper proposed a Discontinuity Capturing Shallow Neural Network (DCSNN)~\cite{HLL22} that allows a single network to represent piecewise smooth functions via a simple augmentation technique. The network is completely shallow (one hidden layer), so the resulting number of trainable parameters is moderate (only a few hundred) and attains prediction accuracy as low as $O(10^{-7})$ in relative $L^2$ norm for all tests in both 2D and 3D elliptic interface problems. Note that the above neural network methods are all completely mesh-free, but their convergence still requires further investigation.

In this work, we propose a novel hybrid method that combines neural network learning machinery and traditional finite difference methods to solve the Poisson interface problem (\ref{Eq:Poisson})-(\ref{Eq:bdc}). The entire computation only comprises a supervised learning task of function approximation and a fast direct solver of the Poisson equation, which can be easily and directly implemented regardless of interface geometry.
\blue{Here, we want to emphasize that it is not our intention to replace traditional numerical methods such as the immersed interface method (IIM) or immersed finite element method (IFEM) nor to compete with them in every aspect. Instead, we want to provide an alternative (especially from the implementation aspect) to solve Poisson interface problems with non-homogeneous jump conditions in which the advantages of using fast Poisson solver and machine learning can be fully exploited. As known, the IIM and non-bodyfitted IFEM need some  complicated treatments to handle the non-homogeneous jump conditions near the interface, especially in 3D case. However, in the present hybrid method, these interface conditions can be easily incorporated into a function constructed by supervised learning and thus regular finite difference scheme can be exploited. The numerical experiments for 2D and 3D Poisson interface problems in Section 3 indicate that the proposed method can achieve a similar accuracy with the IIM.}

The rest of the paper is organized as follows. In Section~2, we present the methodology and list some features, including error analysis of the hybrid scheme. Numerical results for the Poisson interface problems and Stokes equations with singular forces are given in \blue{Sections}~3 and 4, respectively, followed by some concluding remarks and future works in Section~5.

%%%%%%%%%%%%%%%%%%%%%%%%%%%%%%%%%%%%%%%%%
\section{Hybrid neural-network and finite-difference methodology}
%%%%%%%%%%%%%%%%%%%%%%%%%%%%%%%%%%%%%%%%%

By taking advantage of the machine learning techniques, our goal is to design an easy-to-implement fast solver for the Poisson interface problem~(\ref{Eq:Poisson})-(\ref{Eq:bdc}). To this end, we propose a novel hybrid method that exploits the advantages of neural network learning machinery and traditional finite difference method. As we can see from the jump conditions in (\ref{Eq:jump}), the solution $u$ is non-smooth across the interface. Thus, we start by decomposing the solution into
\begin{align}\label{Eq:decomp}
u(\bx) = v(\bx) + w(\bx),
\end{align}
where $v$ and $w$ represent the singular (non-smooth) and regular (smooth) parts of $u$, respectively. More precisely, we require $w$ to be fairly smooth over the entire domain $\Omega$, so that the zero jumps $\dbblk{w} = \dbblk{\partial_n w} = \dbblk{\laplace w} = 0$ on the interface are all satisfied. Now the singular solution $v$ is responsible for having all the discontinuities across the interface; hereby, we construct this discontinuous function by assuming
\begin{align}\label{Eq:v_V}
v(\bx) =
\left\{
\begin{array}{ll}
\mathcal{V}(\bx)  & \bx \in \Omega^-,\\
0          & \bx \in \Omega^+,
\end{array}\right.
\end{align}
where $\mathcal{V}$ is a smooth function to be found. Using the above definition and plugging the decomposition~(\ref{Eq:decomp}) into the jump conditions (\ref{Eq:jump}) and differential equation~(\ref{Eq:Poisson}), the unknown function $\mathcal{V}$ must satisfy the following constraints along the interface:
\begin{align}\label{Eq:V}
\mathcal{V}(\bx) = -\gamma(\bx), \quad \partial_n \mathcal{V}(\bx) = -\rho(\bx), \quad \laplace \mathcal{V}(\bx) = - \dbblk{f(\bx)}, \quad \bx\in\Gamma.
\end{align}
Note that this function is not unique, in the sense that there exist infinitely many functions defined in the domain $\Omega$ that satisfy the restrictions~(\ref{Eq:V}). To find $\mathcal{V}$, we leverage the power of function expressibility of neural networks. Here, we simply employ a shallow (one hidden layer) fully-connected feedforward neural network to approximate $\mathcal{V}$, and learn the function via the supervised learning model. Specifically, given a dataset with $M$ training data points $\{\bx^i_\Gamma\in\Gamma\}_{i=1}^M$ and the target outputs $\gamma(\bx^i_\Gamma)$, $\rho(\bx^i_\Gamma)$ and $\dbblk{f(\bx^i_\Gamma)}$, we find $\mathcal{V}(\bx)$ by minimizing the following mean squared error loss consisting of the
residuals of conditions in Eq.~(\ref{Eq:V}):
\begin{align}\label{Eq:loss}
\Loss(\bp) =\frac{1}{M} \sum_{i=1}^M\left[ \left( \mathcal{V}(\bx_\Gamma^i;\bp)+\gamma(\bx_\Gamma^i) \right)^2 + \left( \partial_n \mathcal{V}(\bx_\Gamma^i;\bp) + \rho(\bx_\Gamma^i) \right)^2 + \left( \laplace \mathcal{V}(\bx_\Gamma^i;\bp) + \dbblk{f(\bx_\Gamma^i)} \right)^2\right],
\end{align}
where $\bp$ collects all trainable parameters (weights and biases) in the network. To train the above loss model, we adopted the Levenberg-Marquardt (LM) method~\cite{Marquardt63}, a full-batch optimization algorithm which is particularly efficient for least squares losses. We should also mention that the partial derivatives of the target function $\mathcal{V}(\bx)$ in the loss function~(\ref{Eq:loss}) can be computed easily by automatic differentiation~\cite{BPRS18}.

Once $\mathcal{V}$ is available, we can obtain $w$ by solving the following Poisson equation:
\begin{align}
&\laplace w(\bx) = \laplace u(\bx) - \laplace v(\bx) =
\left\{
\begin{array}{ll}
f(\bx) - \laplace\mathcal{V}(\bx)  & \bx \in \Omega^-,\\
f(\bx)                                           & \bx \in \Omega^+,
\end{array}\right. \label{Eq:w} \\
&w(\bx) = u_b(\bx), \quad \bx\in\partial\Omega.
\end{align}
Notice that, using the last jump constraint for $\mathcal{V}$ in Eq.~(\ref{Eq:V}), one can immediately see that the right-hand side function of (\ref{Eq:w}) is continuous on the entire domain. Moreover, $w$ is accompanied by exactly the same boundary conditions as the solution $u$, since $v$ vanishes in $\Omega^+$. As a result, $w$ has sufficient regularity (recall the requirement $\dbblk{w} = \dbblk{\partial_n w} = \dbblk{\laplace w} = 0$) and satisfies the Poisson equation in the regular domain $\Omega$ that can be simply and efficiently solved using the well-developed public software Fishpack~\cite{FISHPACK} or any fast Poisson solvers. Of course, other traditional numerical methods, such as finite volume or finite element methods, can also be used to find the solution $w$. Additionally, we remark that both the right-hand sides of Eq.~(\ref{Eq:v_V}) and the Poisson equation~(\ref{Eq:w}) involve the categorization of $\bx$ in $\Omega^-$ or $\Omega^+$, which can be easily done with the assistance of a level set function for which the zero level set represents the interface $\Gamma$.

Let us summarize the proposed hybrid neural-network and finite-difference method for solving the Poisson interface problem~(\ref{Eq:Poisson})-(\ref{Eq:bdc}) as follows:
\paragraph{Step 1} With a given training dataset and a sufficient number of neurons (or trainable parameters) used in the network, find the neural network function $\mathcal{V}$ by minimizing the loss function~(\ref{Eq:loss}) using the LM optimizer. Compute $\mathcal{V}$ at the finite difference discretization grid points in $\Omega^-$ and then obtain the singular part of the solution, $v$, at those grid points using Eq.~(\ref{Eq:v_V}).
\paragraph{Step 2} Evaluate $\laplace\mathcal{V}$ at the grid points in $\Omega^-$, and solve the Poisson equation~(\ref{Eq:w}) by discretizing the Laplace operator using the standard five-point Laplacian and applying a fast direct solver to obtain the regular part of the solution, $w$, at those gird points.
\paragraph{Step 3} Recover the numerical solution $u = v + w$ at the gird points.\\

We conclude this section by introducing several features of the proposed method as follows:
\begin{enumerate}
\item One can immediately deduce that the source of numerical error comes from the network approximation (optimization and network approximation error) and the finite difference approximation (local truncation error). The solution accuracy clearly depends on these two approximations.
\item The right-hand side function of the Poisson equation for $w$ (see Eq.~(\ref{Eq:w})) is continuous, but, in general, has discontinuities in its derivatives across the interface. \red{Under the finite difference discretization in Step~2, we have $\dbblk{w} = \dbblk{\partial_n w} = \dbblk{\laplace w} = 0$ on the interface. So the local truncation error for grid points right adjacent to the interface (irregular points) is $O(h)$ with mesh size $h$, while the one for other grid points (regular points) is $O(h^2)$. Since the number of irregular points is one dimension lower than the number of regular points used in the problem, this $O(h)$ truncation error on an interior set of relative small size does not affect the overall second-order accuracy. (For 1D case, this can be immediately seen by writing the solution with the discrete Green's function (scaled by $h$) so the $O(h)$ local truncation error would make an $O(h^2)$ contribution to the global error.) To prove this, Beale and Layton \cite{BL06} first write this localized $O(h)$ truncation error at irregular points as the discrete divergence of a function which is only $O(h^2)$ in magnitude, and then perform a maximum norm estimate for a discrete elliptic problem with a nonhomogeneous term of divergence form. They are able to prove the solution global error is $O(h^2)$ and its gradient error is
    $O(h^2\log(1/h))$ in maximum norm.
%    can be written as the discrete divergence of a function which is order of $O(h^2)$
%    Since
%these first-order local truncation errors can be regarded as $O(h)$ interior point sources that produce perturbations of magnitude $O(h^2)$ to the numerical solution. As a result, a second-order global error can still be obtained~\cite{LeVeque07}. A more rigorous convergence analysis can be found in the book~\cite{LI2006}.
The present numerical evidence indeed shows that the overall accuracy for solving $w$ is second-order.}
\item With only one-time training, the obtained network function $\mathcal{V}$ is defined in a continuous sense so that it can be used in Step~2 for any grid resolutions.
\item The advantage of traditional grid-based methods is that the boundary conditions are exactly satisfied. The present hybrid method shares the same advantage thanks to the design of Eq.~(\ref{Eq:v_V}). In contrast, most modern deep learning approaches either adopt a penalty term in the loss function (e.g., PINNs~\cite{RPK19} and DCSNN~\cite{HLL22}) or introduce an energy functional to enforce the boundary condition (e.g., shallow Ritz method~\cite{LCLHL22} and deep Nitsche method~\cite{LM21}), leading to an inevitable prediction error along the domain boundary.
\item The proposed hybrid algorithm is easy to implement and efficient. It comprises a supervised learning task (for $\mathcal{V}$) and a fast direct Poisson solver (for $w$), and there are already many well-developed and efficient packages for both tasks. \red{We should also point out that the regular part solver is not limited to finite difference method nor the five-point Laplacian discretization. Finite volume, finite element, or spectral methods can also be used to find the solution $w$.}
\item The present method can be applied to any domain that has available fast solvers, such as a two-dimensional disk, a three-dimensional sphere, higher dimensional cube with periodic boundary conditions.
\item It is straightforward to implement the present method when multiple embedded interfaces are considered.

\end{enumerate}

%%%%%%%%%%%%%%%%%%%%%%%%%%%%%%%%%%%%%%%%%
\section{Numerical results}
%%%%%%%%%%%%%%%%%%%%%%%%%%%%%%%%%%%%%%%%%

In this section, we check the accuracy of the proposed method by performing two numerical tests, including solving two- and three-dimensional Poisson interface problems. In each test, the neural net function $\mathcal{V}$ is simply represented via a shallow network with a sigmoid activation function, in which only a single hidden layer is employed. Thanks to the shallow network structure, it only needs to train a moderate amount of parameters (a few hundred parameters used throughout all numerical examples), so learning this network function is efficient, for example, it can be done in seconds on iMac (2021). Since all the computational domains considered in the following problems are regular (square in 2D and cube in 3D), to solve the regular part $w$, we set up a uniform Cartesian grid layout with the same mesh size $h$ in each spatial direction.

%%%%%%%%%%%%%%%%%%%%%%%%%%%%%%%%%%%%%%%%%
\paragraph{\textbf{Example 1}}
%%%%%%%%%%%%%%%%%%%%%%%%%%%%%%%%%%%%%%%%%

We start by solving a two-dimensional Poisson interface problem and compare the results with the ones obtained by the 2D IIM~\cite{LT08}. The problem is defined in the square domain $\Omega = [-1,1]^2$ in which the embedded interface is an ellipse given by $\Gamma: (x/0.8)^2 + (y/0.2)^2 = 1$. The exact solution is chosen as
\begin{align*}
u_e(x,y) =
\left\{
\begin{array}{ll}
\exp(x)\cos(y)     & \mbox{if\;\;} (x,y) \in \Omega^-,\\
\exp(x^2)\cos(y) & \mbox{if\;\;} (x,y)\in \Omega^+,
\end{array}\right.
\end{align*}
so the corresponding right-hand side $f(x,y)$, and the jump information $\gamma(x,y)$, $\rho(x,y)$ and $\dbblk{f(x,y)}$ used in the loss function can be calculated accordingly. In this example, the network for $\mathcal{V}(x,y)$ is equipped with $40$ neurons in the hidden layer and is trained using $200$ randomly sampled training points on the interface $\Gamma$. We finish the training process when the stopping condition $\Loss(\bp)< 10^{-12}$ or the maximum iteration number (epoch = 1000) is met.

\begin{figure}[t]
\centering
\includegraphics[scale=0.43]{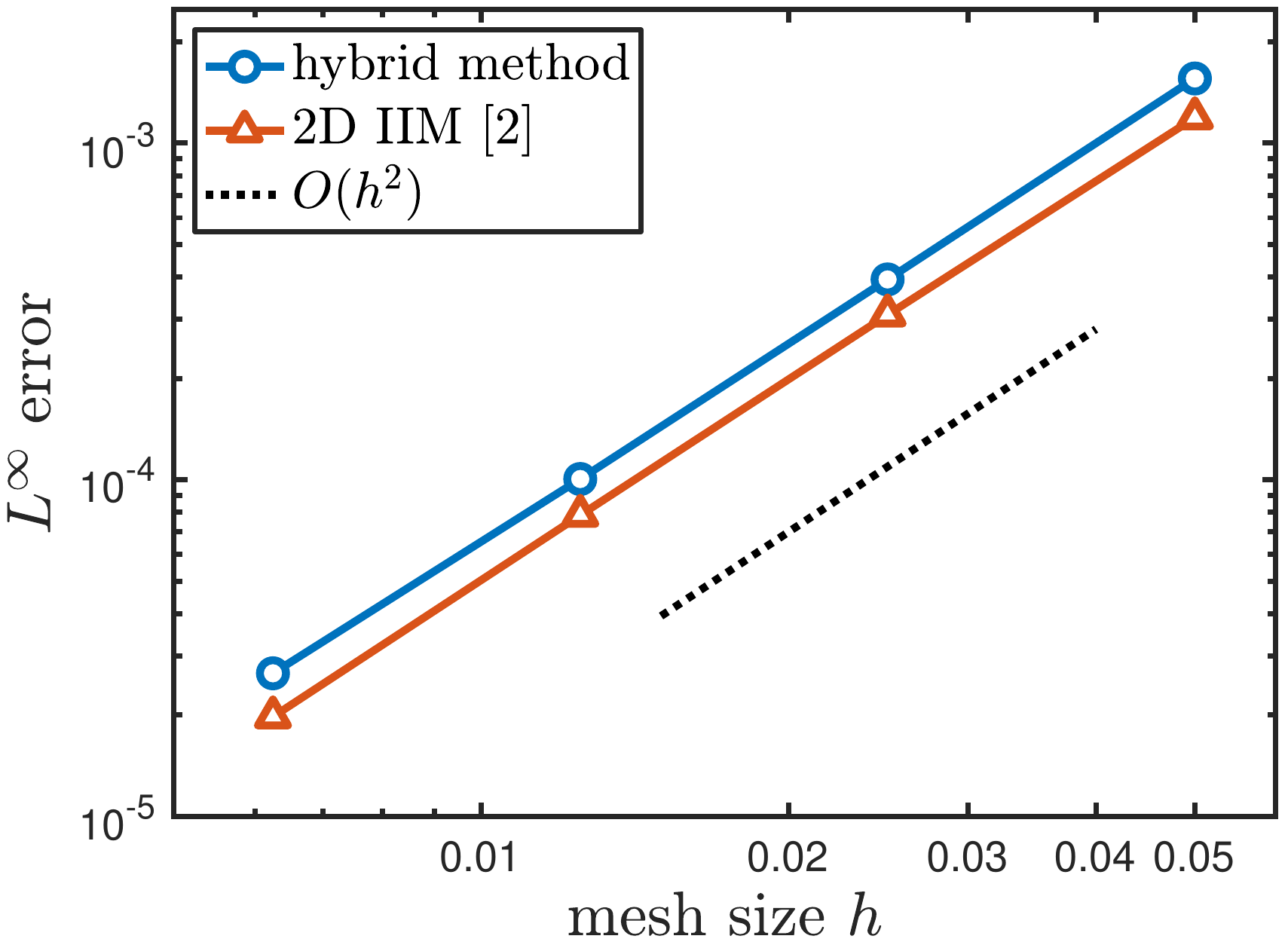}
\includegraphics[scale=0.43]{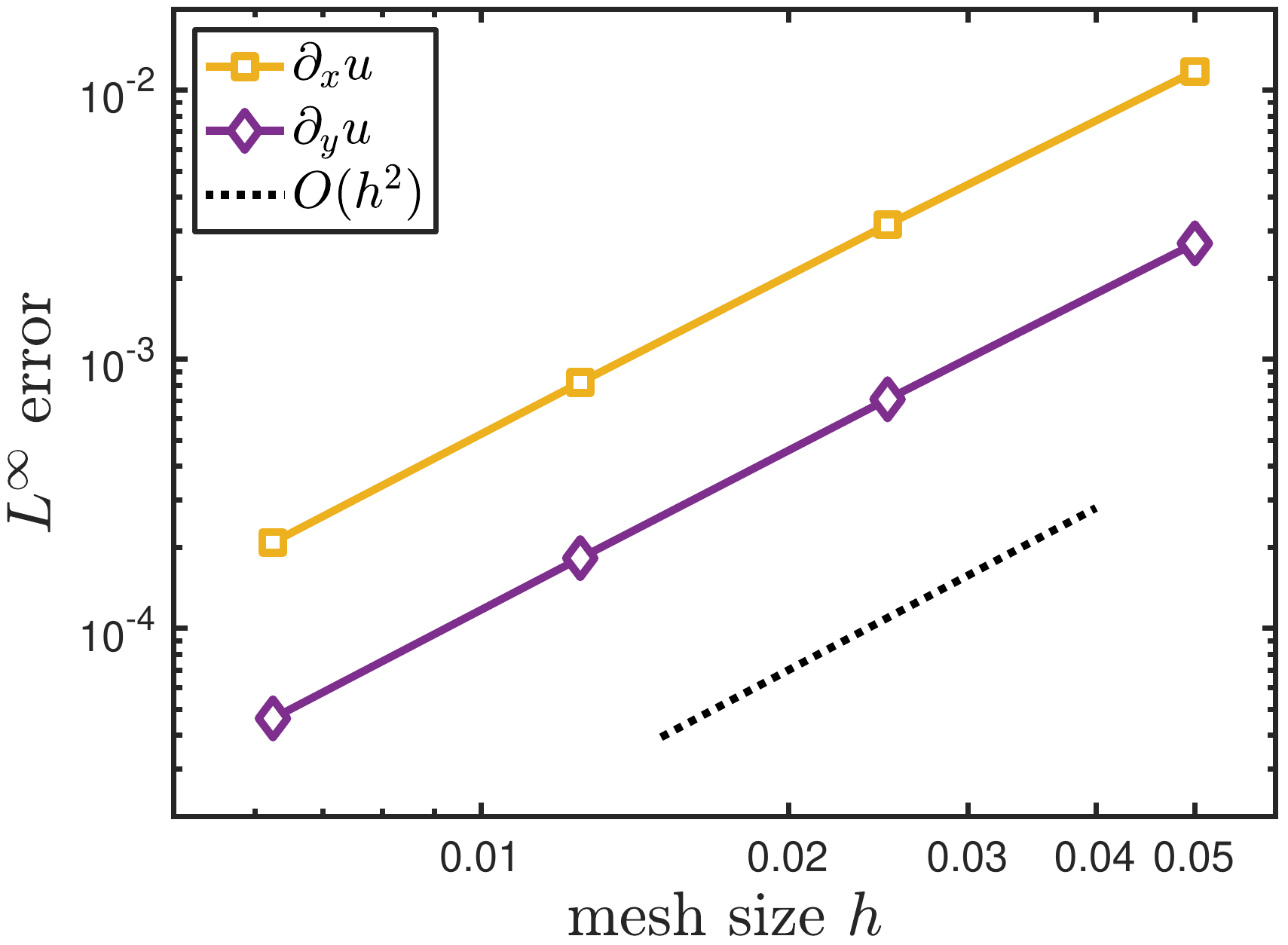}
\caption{Mesh refinement results for the 2D Poisson interface problem in Example 1. Left: Comparison of maximum norm errors of $u$ between the present hybrid method and 2D IIM~\cite{LT08}. Right: Maximum norm error of the gradient $\grad u=(\partial_x u, \partial_y u)$.}
\label{Fig:comparison_2d}
\end{figure}

In the left panel of Fig.~\ref{Fig:comparison_2d}, we report the mesh refinement study for maximum norm error $\|u-u_e\|_\infty$ as a function of mesh size $h$, where $u$ denotes the numerical solution. One can see that the results obtained by the present hybrid method (solid blue line with circular markers) and IIM (solid red line with triangular markers) are almost equally well, and both achieve a second-order convergence rate. We then use the computed solution to find $\grad u=(\partial_x u, \partial_y u)$ simply by applying standard central difference for the regular part $w$ and automatic differentiation for the singular part $v$. As can be seen in the right panel of Fig.~\ref{Fig:comparison_2d}, the gradient of the numerical solution attains a second-order convergence too.

As discussed in the previous section, the induced numerical error comes from both neural network approximation and finite difference approximation. Although not shown here, the final loss value here is about $10^{-13}\sim10^{-14}$, which leads to the predictive accuracy of the target function $\mathcal{V}$ and its Laplacian are of magnitude $10^{-6}\sim10^{-7}$. Thus, we can conclude that the error is mainly dominated by the second-order finite difference approximation error when $h^2 \gtrsim 10^{-6}$.
 To verify our error estimation, we have run more refining numerical tests with $h^2 < 10^{-6}$. As expected, the resulting error cannot be reduced when refining the mesh width $h$ with higher resolutions (not shown here).

%%%%%%%%%%%%%%%%%%%%%%%%%%%%%%%%%%%%%%%%%
\paragraph{\red{\textbf{Example 2}}}
%%%%%%%%%%%%%%%%%%%%%%%%%%%%%%%%%%%%%%%%%
\red{
Because of the mesh-free nature of the neural network, the proposed method is robust to handle complicated interface geometries in which one only has to input the interface description (such as interface parametric form and normal vector) in implementing supervised learning for the singular part solution.
We test the problem where the interface is of super-ellipse shape $(x/\sqrt{0.7})^4 + (y/\sqrt{0.1})^4 = 1$, and choose the same analytic solution $u_e$ and follow the same setups as in Example~1.
In the left panel of Fig.~\ref{Fig:superellipse}, we depict the numerical solution $u$ with the mesh size $h = 1/160$; as one can see that the discontinuity is indeed captured sharply across the interface. In the right panel we report the mesh refinement result for the numerical solution and its gradient. As anticipated, the maximum norm errors converge with second-order accuracy.
}
\begin{figure}[!h]
\centering
\includegraphics[scale=0.43]{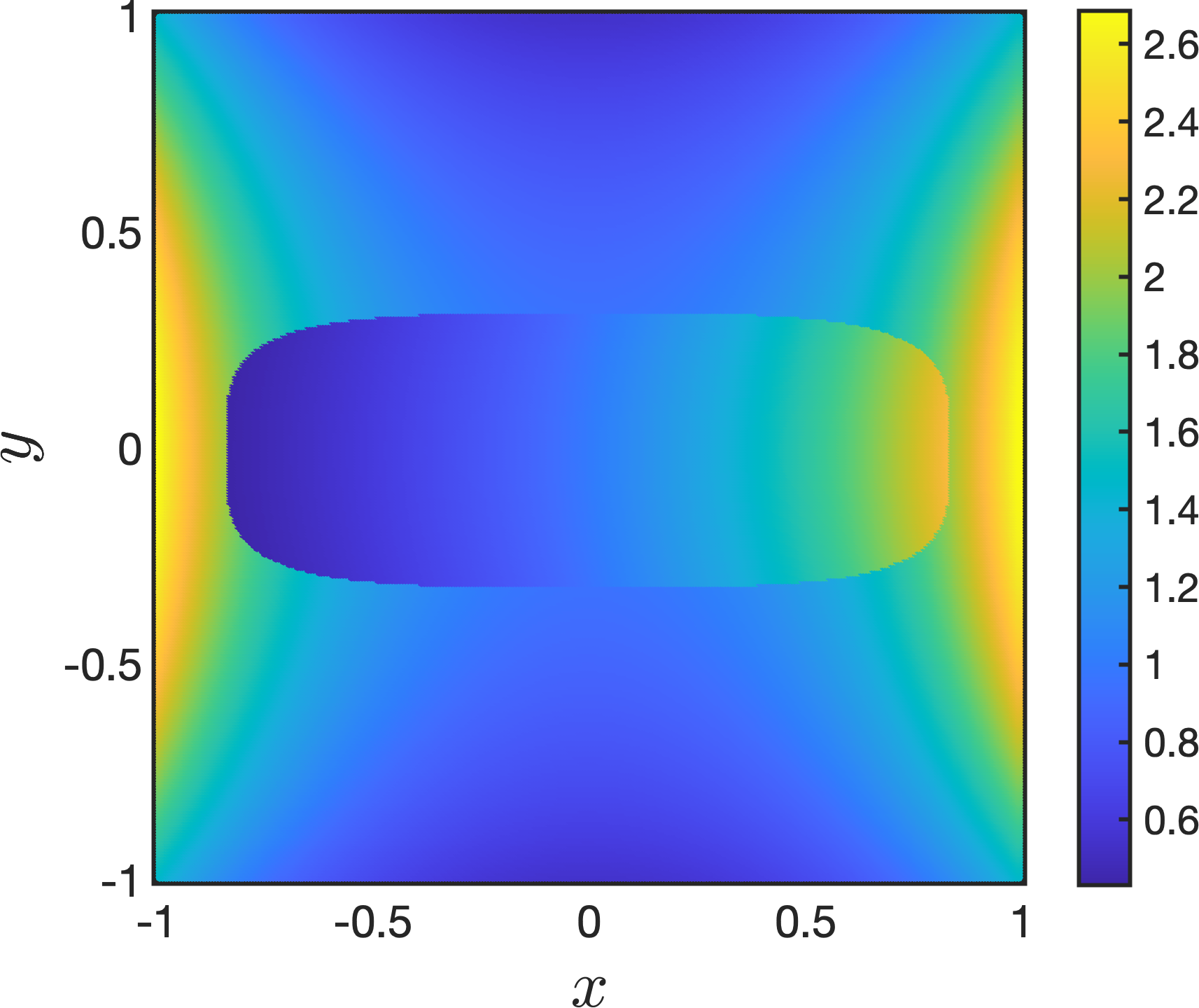}
\includegraphics[scale=0.43]{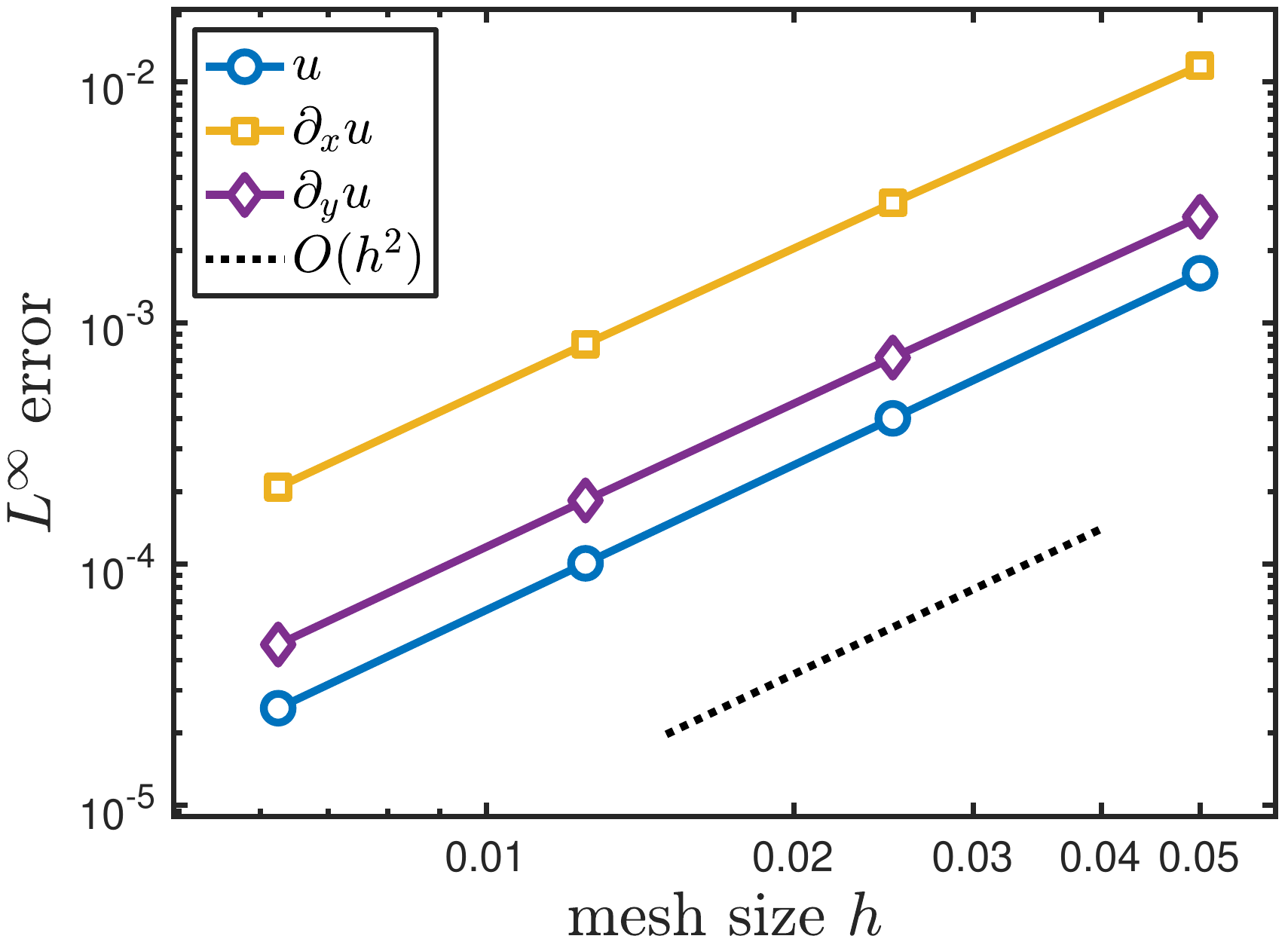}
\caption{Left: The numerical solution $u$ obtained by the present hybrid method with the mesh size $h = 1/160$. Right: Maximum norm errors of the numerical solution $u$ and its gradient $\grad u=(\partial_x u, \partial_y u)$.}
\label{Fig:superellipse}
\end{figure}

%%%%%%%%%%%%%%%%%%%%%%%%%%%%%%%%%%%%%%%%%
\paragraph{\red{\textbf{Example 3}}}
%%%%%%%%%%%%%%%%%%%%%%%%%%%%%%%%%%%%%%%%%
\red{
To showcase the reliability of the hybrid method, we present an example in the scenario that the exact solution is unavailable. We again embed an ellipse $\Gamma$ with parametric form  $(x(s),y(s))=(\sqrt{0.7}\cos s,\sqrt{0.1}\sin s),\,s\in[0,2\pi)$ in the square domain $\Omega = [-1,1]^2$. Specifically, we choose the right-hand side function of the Poisson interface problem as
\begin{align*}
f(x,y) =
\left\{
\begin{array}{ll}
\exp(x\sin y)  & \mbox{ if } (x,y) \in \Omega^-,\\
\exp(y\cos x) & \mbox{ if } (x,y) \in \Omega^+.
\end{array}\right.
\end{align*}
Along the domain boundary $\partial\Omega$, we set the Dirichlet boundary condition $u(x,y) = 0$; along the interface $\Gamma$, we choose the jump condition
$\dbblk{u}(s) = \gamma(s) = \sin s, \dbblk{\partial_n u}(s) = \rho(s) = \cos s$, and the jump $\dbblk{f}(s)$ can be accordingly obtained from the above given function.
The singular part solution is trained using 150 neurons in one hidden layer with 300 randomly sampled training points on the interface $\Gamma$.}

\red{
We need to point out that, since the exact solution is not available in this test, we measure the $L^\infty$ error by the successive error $\|u_h-u_{h/2}\|_\infty$, where $u_h$ denotes the solution with the grid size $h$. We depict the numerical solution $u$ in the left panel of Fig.~3 with the finest resolution $h = 1/320$ in this test. In the right panel we show the maximum norm errors for $u$ and its gradient. Again, all the quantities converge with second-order accuracy.}

%the error of the singular solution is roughly $O(10^{-6})$ (square root of the final loss value), and the regular part can be approximated by the the successive error $\|u_h-u_{h/2}\|_\infty$, where $u_h$ denotes the solution with the grid size $h$. The overall $L^\infty$-error is then the sum of the two errors. We depict the numerical solution $u$ in the left panel of Fig.~3 with the finest resolution $h = 1/320$ in this test. In the right panel we show the maximum norm errors for $u$ and its gradient. Again, all the quantities converge with second-order accuracy.}

\begin{figure}[!h]
\centering
\includegraphics[scale=0.43]{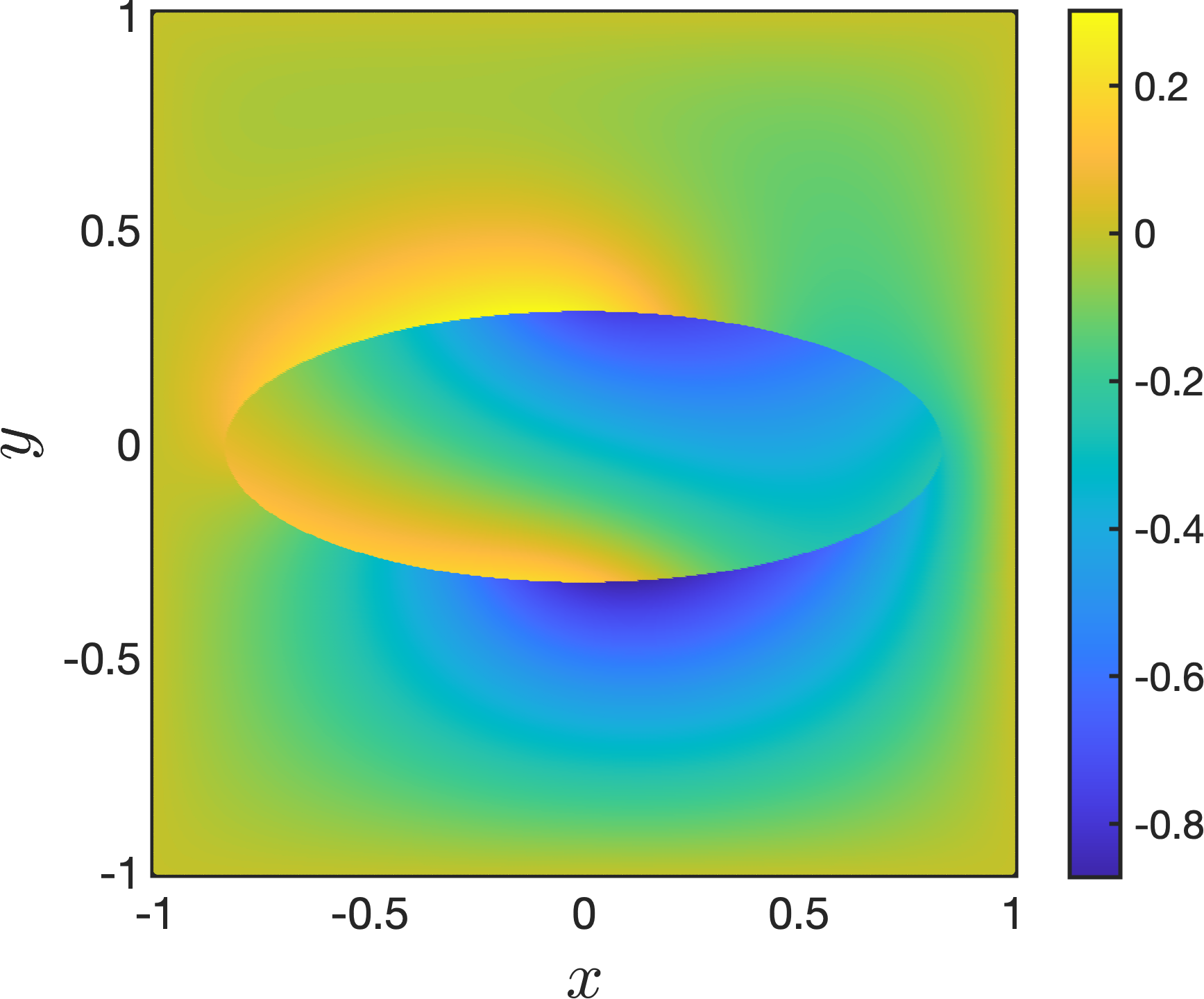}
\includegraphics[scale=0.43]{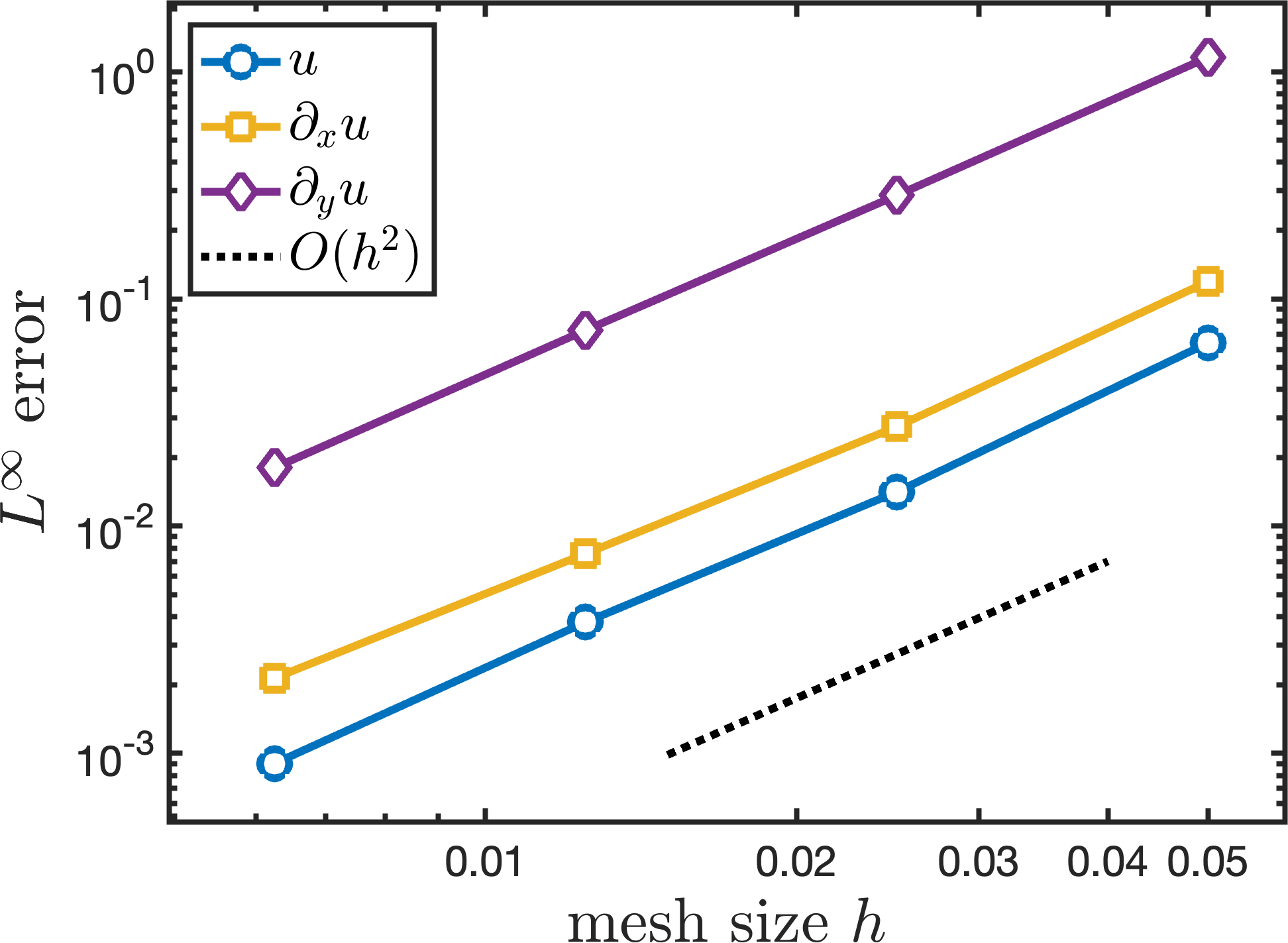}
\setcounter{figure}{2}
\caption{Left: The numerical solution $u$ obtained by the present hybrid method with the mesh size $h = 1/320$. Right: Maximum norm errors of the numerical solution $u$ and its gradient $\nabla u=(\partial_x u, \partial_y u)$.}
\label{Fig:no_exact_solution}
\end{figure}

%%%%%%%%%%%%%%%%%%%%%%%%%%%%%%%%%%%%%%%%%
\paragraph{\textbf{Example 4}}
%%%%%%%%%%%%%%%%%%%%%%%%%%%%%%%%%%%%%%%%%

\begin{figure}[h]
\centering
\includegraphics[scale=0.43]{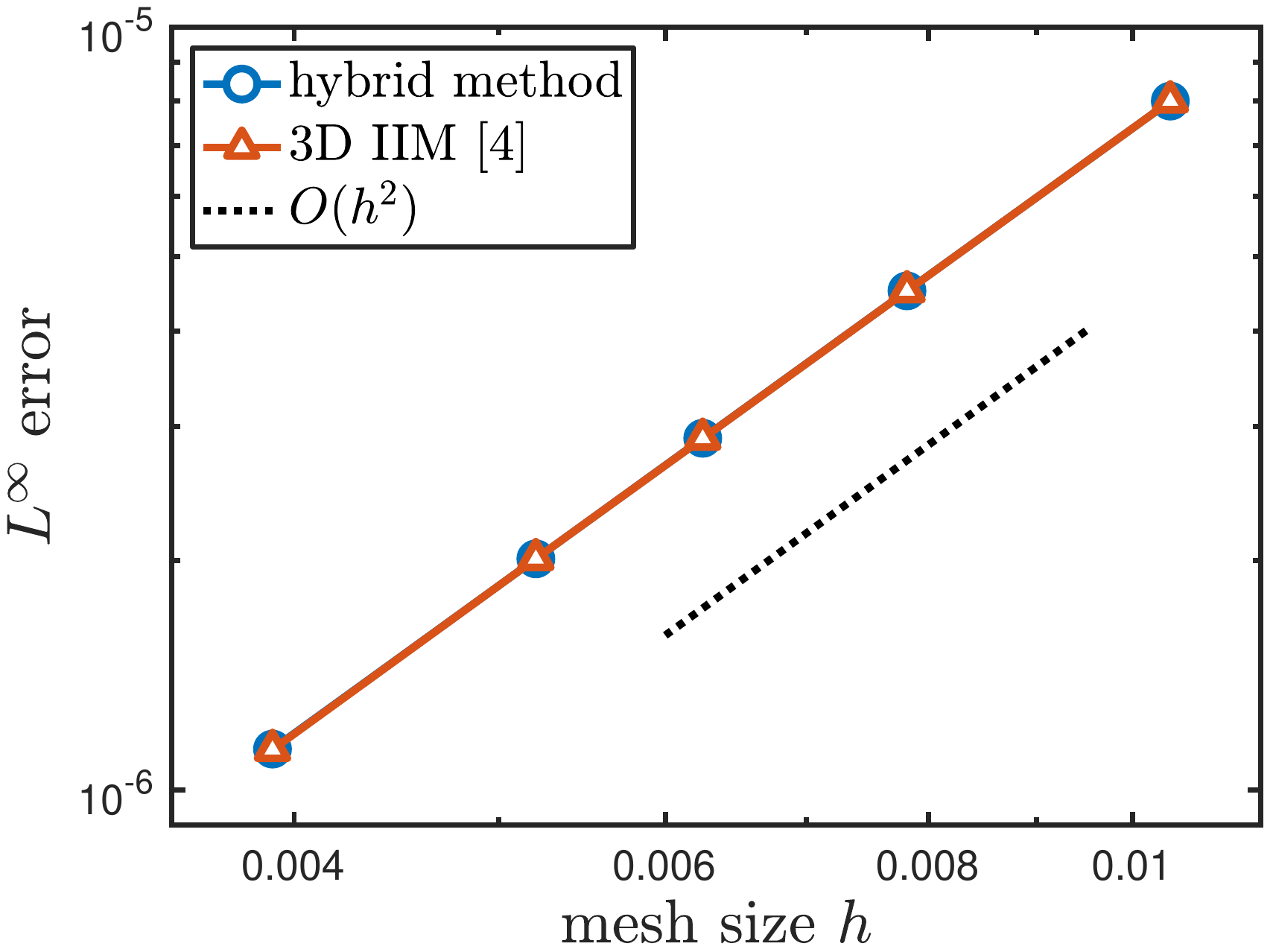}
\includegraphics[scale=0.43]{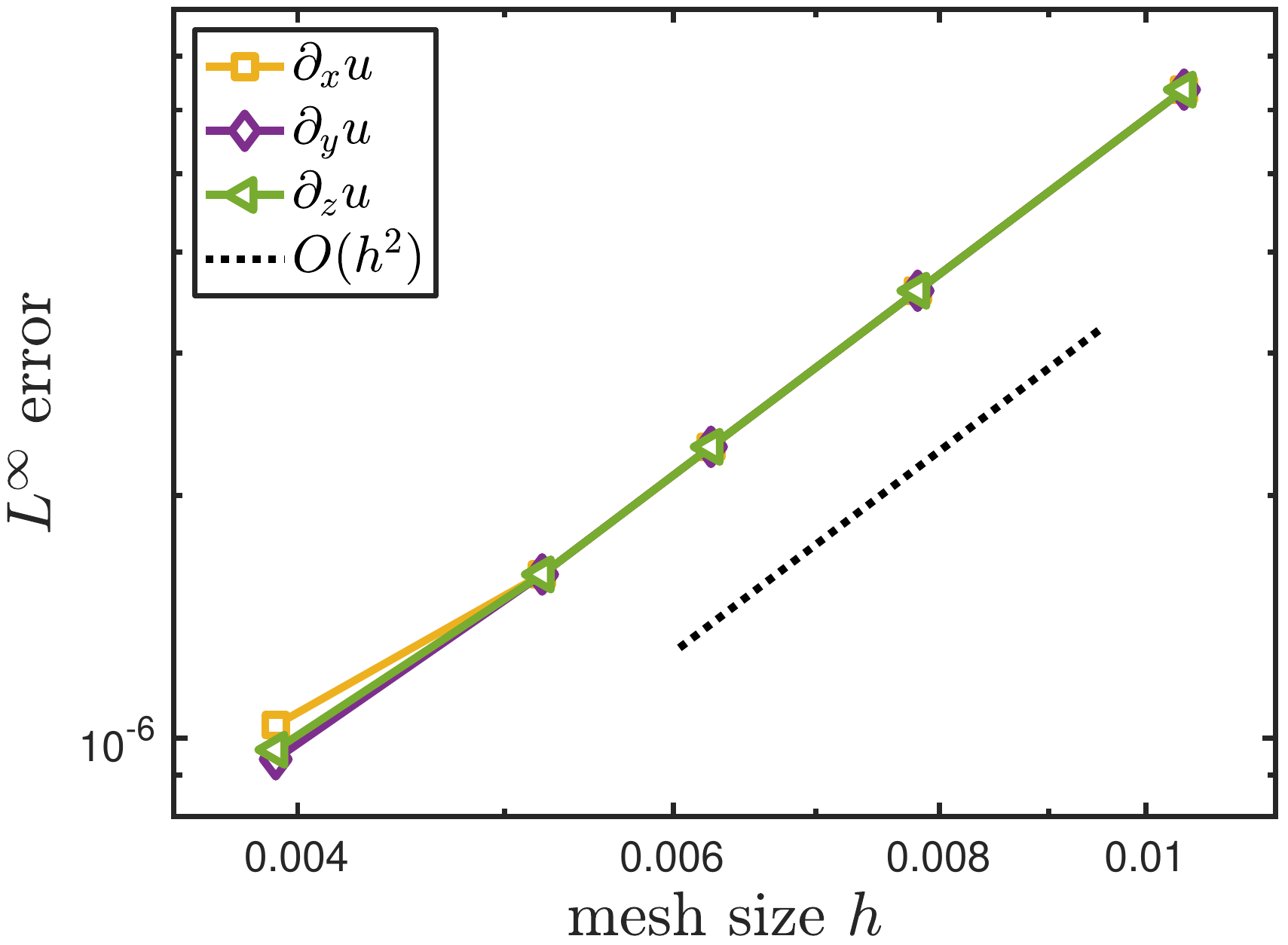}
\caption{Mesh refinement results for the 3D Poisson interface problem in Example 2. Left: Comparison of maximum norm errors of $u$ between the present hybrid method and 3D IIM~\cite{HHL19}. Right: Maximum norm error of the gradient $\grad u = (\partial_x u, \partial_y u, \partial_z u)$.}
\label{Fig:comparison_3d}
\end{figure}

We proceed to consider the three-dimensional Poisson interface problem, in which the interface is an ellipsoid $\Gamma: (x/0.7)^2 + (y/0.5)^2 + (z/0.3)^2 = 1$, embedded in a cube $\Omega = [-1,1]^3$. The exact solution is given by
\begin{align*}
u_e(x,y,z) =
\left\{
\begin{array}{ll}
\exp(x+y+z)            & \mbox{if\;\;} (x,y,z) \in \Omega^-,\\
\sin(x)\sin(y)\sin(z) & \mbox{if\;\;} (x,y,z)\in \Omega^+.
\end{array}\right.
\end{align*}
Again, one can obtain $f$, $\gamma$, $\rho$, and $\dbblk{f}$ accordingly. We use the same shallow network structure as in the previous example, i.e., we set $40$ neurons in one hidden layer and $200$ training points on the interface $\Gamma$ to learn $\mathcal{V}(x,y,z)$. Fig.~\ref{Fig:comparison_3d} shows the mesh refinement results for the present method (solid blue line with circular markers) and 3D IIM solver developed in~\cite{HHL19} (solid red line with triangular markers), as well as the maximum norm errors for the numerical gradient. Similar to the 2D case, one can clearly see that the results of the present hybrid method and IIM are almost identical, and a second-order accuracy is achieved for both the numerical solution and its gradient. It is important to point out that, the implementation for learning $\mathcal{V}$ and solving $w$ in 3D problems is indeed straightforward as in 2D. By contrast, calculating the extra correction terms (incorporating all jump information) in the IIM implementation can be quite tedious in 3D problems, especially when the interface geometry is complex.

%%%%%%%%%%%%%%%%%%%%%%%%%%%%%%%%%%%%%%%%%
\section{Application: 2D Stokes equations with singular forces}
%%%%%%%%%%%%%%%%%%%%%%%%%%%%%%%%%%%%%%%%%
In this section, we apply the proposed method to solve the two-dimensional Stokes equations with singular forces on an interface $\Gamma$, which result from the immersed boundary formulation~\cite{Pes02} of fluid-structure interaction problems. The governing equations are written the same as in \cite{LT08}:
\begin{align}
-\grad p(\bx) + \mu\laplace\bu(\bx) + \int_\Gamma \bF(s) \delta(\bx-\bX(s)) \mbox{ d}s + \bg(\bx) = \mathbf{0}, \quad &\bx\in\Omega, \label{Stokes_1}\\
\div\bu(\bx) = 0, \quad &\bx\in\Omega, \label{Stokes_2}\\
\bu (\bx) = \bu_b(\bx), \quad &\bx\in\partial\Omega, \label{Stokes_3}
\end{align}
where $\bu = (u_1,u_2)$ denotes the fluid velocity, $p$ is the pressure, and $\mu$ is the constant viscosity. There are two different forces acting to the fluid; namely, the external force field $\bg$ (might be discontinuous) to the fluid domain $\Omega$, and the singular force represented by the interfacial force $\bF$ in terms of delta function formulation on  $\Gamma$. Here, the notation  $\bX(s)$ represents the configuration of the interface $\Gamma$ with the arc-length parameter  $s$. (Notice that, the delta function in Eq.~(\ref{Stokes_1}) is two-dimensional but the integration is over the one-dimensional interface which leaves the integral term has one-dimensional singularity.)  Furthermore, the interfacial force $\bF(s)$ can be written as a sum of the tangential ($\btau$) and normal ($\bn$) components as $\bF = F_\tau\btau + F_n\bn$. Instead of using the delta function formulation, one can rewrite Eqs.~(\ref{Stokes_1})-(\ref{Stokes_3}) into the following Stokes equations with jump conditions across the interface (see the reference in \cite{LL97}):
\begin{align}
-\grad p(\bx) + \mu\laplace\bu(\bx)  + \bg(\bx) = \mathbf{0}, \quad
&\bx\in \Omega^- \cup \Omega^+, \label{Stokes_11}\\
\div\bu(\bx) = 0, \quad &\bx\in\Omega^- \cup \Omega^+, \label{Stokes_21}\\
\bu (\bx) = \bu_b(\bx), \quad &\bx\in\partial\Omega, \label{Stokes_31}\\
\dbblk{p(\bX(s))}=F_n(s), \, \dbblk{\bu(\bX(s))}=\mathbf{0}, \, \dbblk{\partial_n\bu(\bX(s))}= -F_\tau(s)\btau(s)/\mu, \quad
&\bX(s) \in \Gamma\label{Stokes_41}.
\end{align}
To solve the above Stokes equations by the present methodology given in Section 2, following the implementation in~\cite{LL97}, we
firstly apply the divergence operator to Eq.~(\ref{Stokes_11}) and use the divergence-free condition (\ref{Stokes_21}) to obtain the pressure Poisson equation as
\begin{align}
\laplace p (\bx)= \div\bg(\bx), \quad &\bx\in \Omega^- \cup \Omega^+, \label{pressure}\\
\dbblk{p(\bX(s))} = F_n(s), \quad \dbblk{\partial_n p(\bX(s))} = \frac{\pd F_\tau}{\pd s}+ \dbblk{\bg(\bX(s))} \cdot\bn(s), \quad &\bX(s)\in\Gamma.
\label{pressure_jump}
\end{align}
One can immediately see that an extra normal derivative jump condition of the pressure is needed since now the equation (\ref{pressure}) involves second-order derivatives rather than the first-order derivatives in Eq.~(\ref{Stokes_11}). For completeness, we list the derivation of the jump conditions in Eq.~(\ref{pressure_jump}) in Appendix, in which the proof is slightly different from the one given in \cite{LL97}. Once the pressure is found, the velocity can be obtained by solving
\begin{align}
\laplace \bu (\bx) = \frac{1}{\mu}(\grad p(\bx) -\bg(\bx)), \quad &\bx\in \Omega^- \cup \Omega^+,\label{eq:vel_1}\\
\dbblk{\bu(\bX(s))} = \mathbf{0}, \quad \dbblk{\partial_n\bu(\bX(s))} = -F_\tau(s)\btau(s)/\mu, \quad & \bX(s)\in\Gamma, \\
\bu(\bx) = \bu_b(\bx), \quad &\bx\in\partial\Omega.
\end{align}
One can immediately see that these are Poisson interface problems so the present method can be applied directly.

\paragraph{\textbf{Example}}
For the numerical test, we use the example as in~\cite{LT08}. We consider a square domain $\Omega = [-2,2]^2$ and the interface is simply a unit circle with the center located at the origin: $\Gamma = \{\bX(s) = (\cos s,\sin s)|s \in[0,2\pi)\}$. The interfacial force comprises both tangential ($\btau(s)=\bX'(s)$) and normal ($\bn(s)= \bX(s)$) directions as $\bF(s) = 2\sin(3s) \btau(s)- \cos^3(s)\bn(s)$. Here, we set $\mu = 1$, and the exact velocity written in polar coordinates is chosen as
\begin{align*}
&u_1(r,\theta) =
\left\{
\begin{array}{ll}
\frac{1}{8}r^2\cos(2\theta)+\frac{1}{16}r^4\cos(4\theta)-\frac{1}{4}r^4\cos(2\theta)             & \mbox{if\;\;} r<1,\\[.2cm]
-\frac{1}{8}r^{-2}\cos(2\theta)+\frac{5}{16}r^{-4}\cos(4\theta)-\frac{1}{4}r^{-2}\cos(4\theta) & \mbox{if\;\;} r\geq1,
\end{array}\right.\\[.2cm]
&u_2(r,\theta) =
\left\{
\begin{array}{ll}
-\frac{1}{8}r^2\sin(2\theta)+\frac{1}{16}r^4\sin(4\theta)+\frac{1}{4}r^4\sin(2\theta)          & \mbox{if\;\;} r<1,\\[.2cm]
\frac{1}{8}r^{-2}\sin(2\theta)+\frac{5}{16}r^{-4}\sin(4\theta)-\frac{1}{4}r^{-2}\sin(4\theta) & \mbox{if\;\;} r\geq1,
\end{array}\right.
\end{align*}
while the exact pressure written in Cartesian coordinates $(x = r\cos(\theta), y = r\sin(\theta))$ is
\begin{align*}
p(x,y) =
\left\{
\begin{array}{ll}
x^3 + \cos(\pi x)\cos(\pi y) & \mbox{if\;\;} r<1,\\
\cos(\pi x)\cos(\pi y)           & \mbox{if\;\;} r\geq1.
\end{array}\right.
\end{align*}
Substituting the above velocity $\bu=(u_1,u_2)$ and pressure $p$ into Eq.~(\ref{Stokes_11}), one can obtain the corresponding discontinuous external force $\bg$.

We establish a three-ouput neural network, $\boldsymbol{\mathcal{V}} = (\mathcal{V}_{u_1}, \mathcal{V}_{u_2}, \mathcal{V}_p)$, to learn the target functions simultaneously. Consequently, the loss function now consists all the jump residuals for $\bu$ and $p$. Notice that, one requires the jump information of the right-hand side of Eq.~(\ref{eq:vel_1}), $\frac{1}{\mu}(\dbblk{\grad p}-\dbblk{\bg})$; the latter is given directly and the former can be obtained by using the identity $\dbblk{\grad p} = \pd_s\dbblk{p} \btau + \dbblk{\partial_n p}\bn$ that links to the given quantities using Eq.~(\ref{pressure_jump}). In the present test, we use a one-hidden-layer network that consists $50$ neurons in the hidden layer. We use $200$ randomly distributed sampling points and train the model using the LM optimizer.

Once the singular parts are obtained, we use fast Poisson solver to find the regular parts.
Inside $\Omega$, the fluid variables are defined at usual staggered grid layout with uniform mesh width $h$. Namely, the velocity components $u_1$ and $u_2$ are correspondingly defined at the cell edges
\[
(x_{i-1/2},y_j) = (-2+(i-1)h,-2+(j-1/2)h), \quad (x_i,y_{j-1/2}) = (-2+(i-1/2)h,-2+(j-1)h),
\]
while the pressure $p$ is defined at the cell center
\[
(x_{i-1/2},y_{j-1/2}) = (-2+(i-1/2)h,-2+(j-1/2)h).
\]
We should point out that, for testing purpose, the pressure boundary condition is chosen to be the Neumann type (which is commonly used in projection method for the pressure increment in Navier-Stokes flows~\cite{BCM01}). After the numerical solution for the regular part of pressure is obtained, we compute  $\grad p$ at the cell edges (which coincides with the location of $\bu$) by applying standard central difference for the regular part and auto differentiation for the singular part.

\begin{figure}[h]
\centering
\includegraphics[scale=0.42]{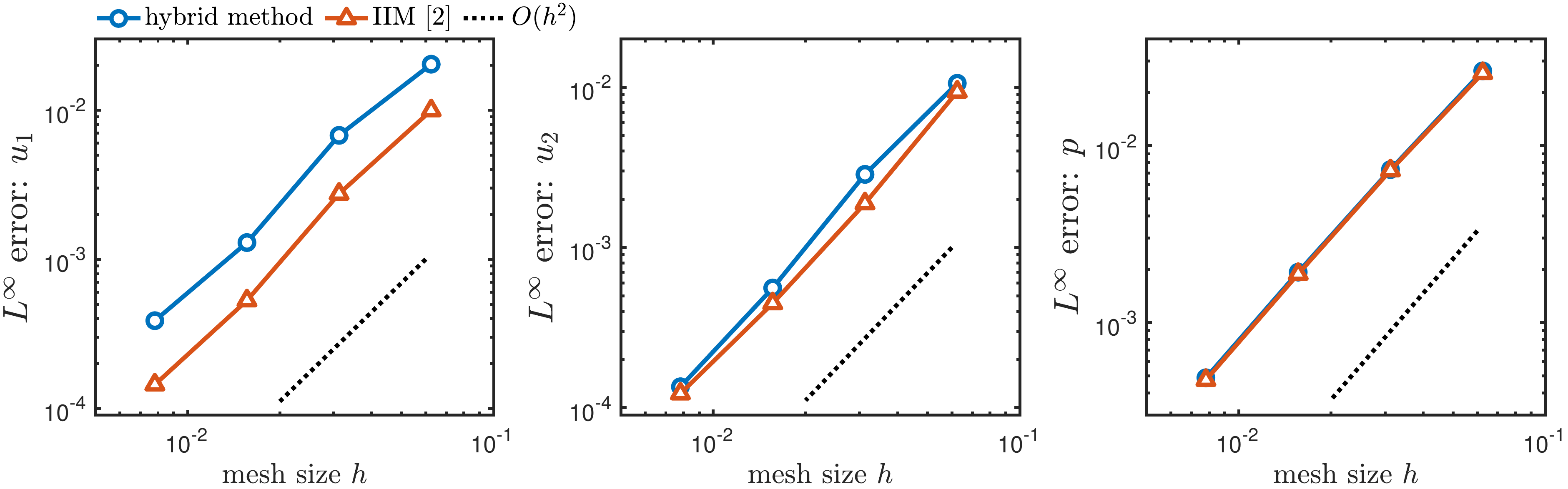}\\[4mm]
\includegraphics[scale=0.42]{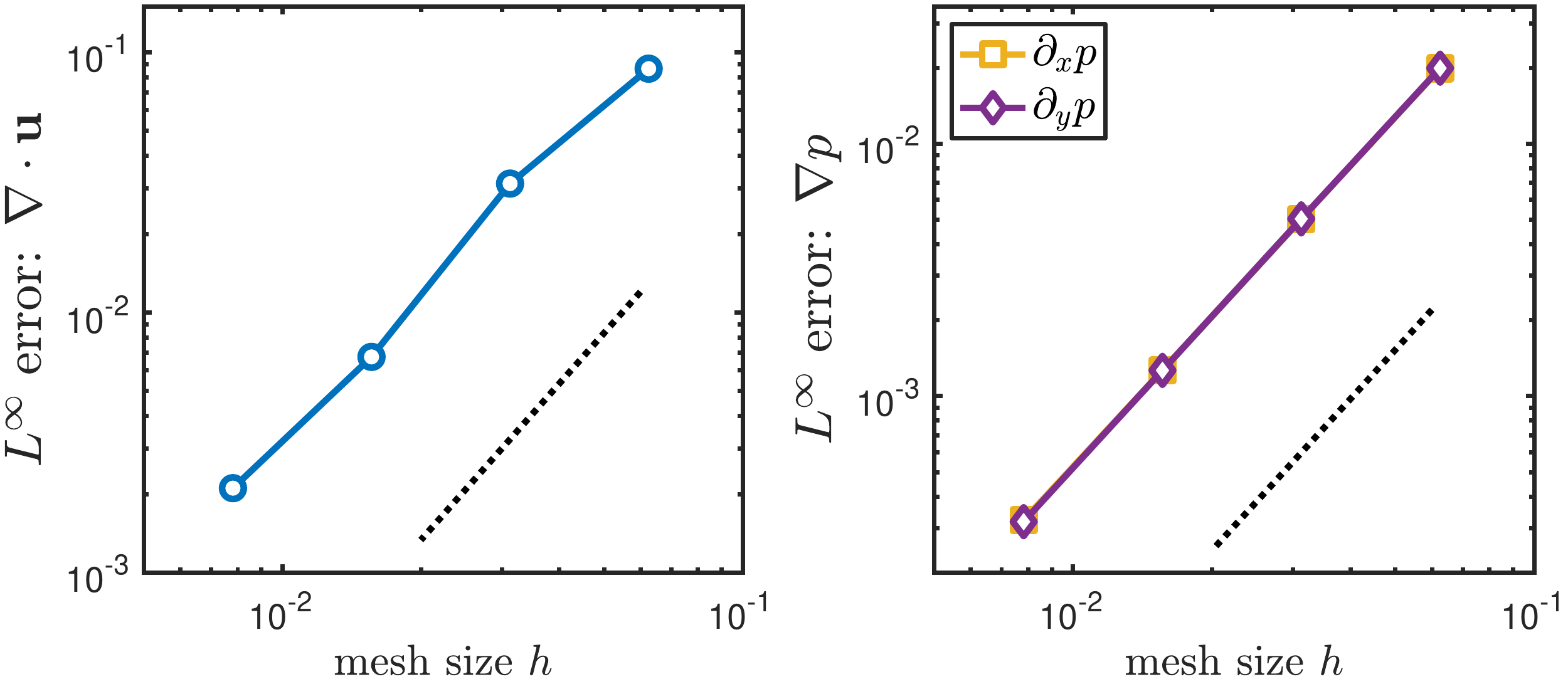}
\caption{Mesh refinement results for 2D Stokes equations with singular forces.  Top: Comparison results between the present hybrid method and 2D IIM~\cite{LT08}. Bottom: Maximum norm errors of $\div\bu$ and $\grad p$.}
\label{Fig:Stokes}
\end{figure}

The resulting errors for all fluid variables are shown in the top panel of Fig.~\ref{Fig:Stokes}. As expected, the accuracy obtained by the present hybrid method is comparable with the ones by IIM~\cite{LT08}, and roughly achieves second-order convergence in maximum norm error. In the bottom panel, we can also see that the numerical errors of $\div\bu$ and $\grad p$ are around second-order as well.
%It should be noticed that the term $\div\bu$ would not vanish numerically since the present solver for Stokes equations only guarantees $\laplace(\div\bu) = 0$ numerically.

%%%%%%%%%%%%%%%%%%%%%%%%%%%%%%%%%%%%%%%%%
\section{Conclusion and future work}
%%%%%%%%%%%%%%%%%%%%%%%%%%%%%%%%%%%%%%%%%

In this paper, we propose a new class of numerical methods to solve an elliptic interface problem whose solution and derivatives are known to have jump discontinuities across an interface. The crucial idea is to decompose the solution into singular (non-smooth) and regular (smooth) parts. The singular part is formed by a neural network representation found by using supervised learning machinery that incorporates all given jump information into the loss function. The regular part, however, is a solution to the Poisson equation, which can be obtained efficiently by several well-developed numerical methods, such as the fast direct solver based on finite difference discretization. Therefore, it is simple to implement our proposed method, and it is straightforward to handle multiple interfaces or high-dimensional problems. The numerical experiments for 2D and 3D Poisson interface problems show that the proposed neural-network and finite-difference hybrid method can achieve second-order accuracy for the solution and its derivatives. Although all illustrated examples consider a single embedded interface only, it is straightforward to implement the hybrid method with multiple interfaces. As an application, we use the present methodology to solve 2D Stokes equations with singular forces. Again, the numerical result shows that all the fluid variables and their derivatives have second-order convergence in maximum norm error as well.

The present hybrid method readily serves as a fast solver for Poisson interface problems involved in the projection step of Navier-Stokes flow problems. \blue{Our future work aims to extend the present methodology for solving variable-coefficient elliptic interface problems in regular or even irregular domains.}

%%%%%%%%%%%%%%%%%%%%%%%%%%%%%%%%%%%%%%%%%
\section*{Acknowledgement}
%%%%%%%%%%%%%%%%%%%%%%%%%%%%%%%%%%%%%%%%%

W.-F. Hu, T.-S. Lin, Y.-H. Tseng, and M.-C. Lai acknowledge the supports by National Science and Technology Council, Taiwan, under the research grants 111-2115-M-008-009-MY3, 111-2628-M-A49-008-MY4, 111-2115-M-390-002, and 110-2115-M-A49-011-MY3, respectively. W.-F. Hu and T.-S. Lin also acknowledge the supports by National Center for Theoretical Sciences,  Taiwan.

\section*{Appendix}
In this appendix, we present the derivation of the pressure jump condition in Eq.~(\ref{pressure_jump}). Let us recall that the domain  $\Omega \subset \mathbb{R}^d$, is separated by an embedded interface $\Gamma$ such that the subdomains inside and outside the interface are denoted by $\Omega^-$ and $\Omega^+$, respectively. And the shorthand $\partial_n$ represents the normal derivative of a quantity where $\bn$ is the unit outward normal vector pointing from $\Omega^-$ to $\Omega^+$. By taking the divergence operator  to Eq.~(\ref{Stokes_1}) and using the divergence-free condition (\ref{Stokes_2}), we obtain the pressure Poisson equation as
\[
\Delta p(\bx) = \div \int_\Gamma \bF(s) \delta(\bx-\bX(s)) \mbox{ d}s + \div \bg(\bx). \]
Since the right-hand side of the above equation involves taking the divergence operator on the Dirac delta function and discontinuous function $\bg$, we should regard them in the sense of distributions. In other words, the above equation should be represented as $\langle \Delta p, \phi \rangle =\langle p, \Delta \phi \rangle$, for all test functions $\phi \in C^{\infty}_0(\Omega)$. So in the following derivations, the test function $\phi$ and its normal derivative $\pd_n \phi$ will be vanished on the outside boundary $\pd \Omega^+$. Applying the derivative properties of the Dirac delta function, we have
\beqs
\langle \Delta p, \phi \rangle & = & \int_\Omega \div \int_\Gamma \bF(s) \delta(\bx-\bX(s)) \mbox{ d}s \,\phi(\bx) \mbox{ d}\bx +
\langle \div \bg, \phi \rangle \nonumber \\
% \int_\Omega \div \bg(\bx) \, d\bx\\
& = & - \int_\Gamma \bF(s) \cdot \grad \phi(\bX(s)) \mbox{ d}s - \langle \bg, \grad \phi \rangle \quad (\mbox{by definition of derivative distributions}) \nonumber \\
&=&  - \int_\Gamma (F_\tau\btau + F_n\bn) \cdot \grad \phi \mbox{ d}s - \int_{\Omega^+} \bg \cdot \grad \phi \mbox{ d}\bx - \int_{\Omega^-} \bg \cdot \grad \phi \mbox{ d}\bx  \nonumber \\
&=& - \int_\Gamma (F_\tau\btau + F_n\bn) \cdot \grad \phi \mbox{ d}s + \int_{\Omega^+} (\div \bg)\phi \mbox{ d}\bx
+
\int_\Gamma (\bg^+ \cdot \bn) \phi \mbox{ d}s \nonumber \\
& & +  \int_{\Omega^-} (\div \bg)\phi \mbox{ d}\bx -
\int_\Gamma (\bg^- \cdot \bn) \phi \mbox{ d}s \quad (\mbox{by applying Green's first identity to $\Omega^\pm$ separately}) \nonumber \\
%& =& - \int_\Gamma (F_\tau \pd_\btau \phi+ F_n \pd_n \phi) \, ds + \int_{\Omega} \div \bg \phi d\bx + \int_\Gamma \dbblk{\bg} \cdot \bn \, ds \\
&=& \int_\Gamma \frac{\pd F_\tau}{\pd s}\phi \mbox{ d}s - \int_\Gamma F_n \,\pd_n\phi \mbox{ d}s + \int_{\Omega^- \cup \Omega^+} (\div \bg) \phi \mbox{ d}\bx + \int_\Gamma \left(\dbblk{\bg} \cdot \bn\right)\phi \mbox{ d}s, \label{LHS}
\eeqs
where the first term in the last equation is obtained using integration by parts.
Meanwhile,
\beqs
\langle p, \Delta \phi \rangle & =& \int_{\Omega^+} p\, \Delta \phi \mbox{ d}\bx + \int_{\Omega^-} p \,\Delta \phi \mbox{ d}\bx \nonumber\\
&=& - \int_{\Omega^+} \grad p \cdot \grad \phi \mbox{ d}\bx - \int_\Gamma p^+ \pd_n \phi \mbox{ d}s
- \int_{\Omega^-} \grad p \cdot \grad \phi \mbox{ d}\bx + \int_\Gamma p^- \pd_n \phi \mbox{ d}s \nonumber \\
&=& - \int_{\Omega^+} \grad p \cdot \grad \phi \mbox{ d}\bx - \int_{\Omega^-} \grad p \cdot \grad \phi\mbox{ d}\bx - \int_\Gamma \dbblk{p}\,
\pd_n \phi \mbox{ d}s \nonumber\\
&=& \int_{\Omega^+} \Delta p \, \phi \mbox{ d}\bx  + \int_{\Gamma} \pd_n p^+  \phi \mbox{ d}s
+ \int_{\Omega^-} \Delta p \, \phi \mbox{ d}\bx - \int_{\Gamma} \pd_n p^-\phi \mbox{ d}s - \int_\Gamma \dbblk{p}\,
\pd_n \phi \mbox{ d}s \nonumber \\
&=& \int_{\Omega^+} (\div \bg) \phi \mbox{ d}\bx + \int_{\Omega^-} (\div \bg) \phi \mbox{ d}\bx
%&=& \int_{\Omega} \Delta p \, \phi \,d\bx
+ \int_\Gamma \dbblk{\pd_n p}
 \phi \mbox{ d}s - \int_\Gamma \dbblk{p}\pd_n \phi \mbox{ d}s \nonumber \\
&=& \int_{\Omega^- \cup \Omega^+} (\div \bg) \phi \mbox{ d}\bx + \int_\Gamma \dbblk{\pd_n p}
 \phi \mbox{ d}s - \int_\Gamma \dbblk{p}\,\pd_n \phi \mbox{ d}s. \label{RHS}
\eeqs
By equating Eq.~(\ref{LHS}) and (\ref{RHS}), one can immediately obtain the jump conditions
\[
\dbblk{p} = F_n, \qquad \dbblk{\pd_n p}= \frac{\pd F_\tau}{\pd s} + \dbblk{\bg} \cdot \bn. \]

%%%%%%%%%%%%%%%%%%%%%%%%%%%%%%%%%%%%%%%%%
%\section*{References}
%%%%%%%%%%%%%%%%%%%%%%%%%%%%%%%%%%%%%%%%%

%\bibliography{hybrid.bib}

\end{document}